\documentclass{article}
\usepackage{amsmath}
\newtheorem{theorem}{Theorem}[section]
\newtheorem{lemma}[theorem]{Lemma}
\newtheorem{proposition}[theorem]{Proposition}
\newtheorem{corollary}[theorem]{Corollary}
\newtheorem{conjecture}[theorem]{Conjecture}

\newenvironment{claim}[1][\it{Claim}]{\begin{trivlist}
\item[\hskip \labelsep {\bfseries #1}]}{\end{trivlist}}
\newenvironment{proof}[1][\it{Proof.}]{\begin{trivlist}
\item[\hskip \labelsep {\bfseries #1}]}{\end{trivlist}}

\newenvironment{remark}[1][\it Remark]{\begin{trivlist}
\item[\hskip \labelsep {\bfseries #1}]}{\end{trivlist}}
\newcommand{\qed}{\nobreak \ifvmode \relax \else
      \ifdim\lastskip<1.5em \hskip-\lastskip\
      \hskip1.5em plus0em minus0.5em \fi \nobreak
      \vrule height0.75em width0.5em depth0.25em\fi}

\title{ Uniform  Annihilators  of Local Cohomology }
\author{{  Caijun  Zhou    }
\\
{\small {\it Department of  Mathematics, Shanghai Normal University,}}\\
{\small {\it Shanghai, 200234, China}}}
\date{}

\begin{document}
\maketitle

\section {Introduction}

 Throughout this paper all  rings  are commutative, associative with
identity, and noetherian. For any unexplained notation and
terminology  we refer the reader to [Ma].

Recall that  a local noetherian ring $R$ is equidimensional if
$\text{dim}(R)  =  \text{dim}(R/P)$ for all minimal primes $P$ of
$R$, and thus a noetherian  ring $R$  is said to be locally
equidimensional if  $R_m$ is   equidimensional for every maximal
ideal $m$ of $R$.  We will use $R^{\circ}$ to denote the
complement of the union of the minimal primes of $R$.

Let $R$ be a   noetherian  ring and $I$ be  an  ideal $I$. For
a $R$-module $M$, we will write   $\text{H}^i_I(M)$ for the $i$-th
local cohomology module of $M$  with support in $\text{V}(I) =
\{P\in \text{Spec}(R)\mid P\supseteq I\}$. We will say an element
$x$ of $R^{\circ}$ is a  uniform local cohomological annihilator of
$R$, if for every maximal ideal $m$, $x$ kills $\text{H}_m^i(R)$ for
all $i$ less than the height of $m$. Moreover, we say that $x$ is a
strong uniform local cohomological annihilator  of $R$ if $x$ is a
uniform local cohomological annihilator  of $R_P$ for every non
minimal prime ideal $P$ of $R$.

It is well known that  a nonzero   local cohomology   module   is
rarely finitely generated, even  the   annihilators of it are not
known in general. On the other hand, it has been discovered by
Hochster and Huneke  that the existence of  a  uniform  local
cohomological annihilator is of great importance in  solving  the
problems such as the existence big Cohen-Macaulay algebras [HH2]
and a uniform  Artin-Rees theorem [Hu]. So it is very interesting
to find out more  noetherian rings  containing uniform local
cohomological annihilators.

A traditional way of studying uniform local cohomological
annihilators is to make use of the dualizing complex over a local
ring. Roberts  initiated  this method in [Ro]. By means of this
technique Hochster and Huneke [HH1] prove that if a locally
equidimensional   noetherian ring $R$ is a homomorphic image
of a Gorenstein ring of finite dimension, then $R$ has a strong
uniform local cohomological   annihilator. It is known not every
local ring has a dualizing complex, however, by passing to
completion, Hochster and   Huneke [HH2] show that an unmixed,
equidimensional  excellent local ring has a strong uniform
local cohomological annihilator.

It is worth noting that a lot of results concerning the
annihilators of local cohomology modules have been established in
recent years. Let us recall some notions before we state these
achievements. One can refer to [BS] for details.

Given ideals $I, J$ in a noetherian ring $R$ and a $R$-modules
$M$, we set

$$\lambda_I^J(M)=\text{inf}\{\text{depth}(M_P)+ \text{ht}(I+P/P)\mid
P\in \text{Spec}(R)\setminus \text{V}(J)\}$$

$$f_I^J(M)=\text{inf}\{i\mid J^n\text{H}_I^i(M)\neq 0 \ \ \text{for \  all \ positive  \ integers  }
n\}.$$ We will say that the Annihilator Theorem for local
cohomology modules holds over $R$ if $\lambda_I^J(M)=f_I^J(M)$ for
every choice of the finitely generated $R$-module $M$ and for
every choice of ideals $I, J$ of $R$. We also say that the
Local-global Principle for the annihilation of local cohomology
modules holds over $R$ if
$f_I^J(M)=\text{inf}\{f_{IR_P}^{JR_P}(M_P)\mid P\in
\text{Spec}(R)\}$ holds for every choice of ideals $I, J$ of $R$
and for every choice of the finitely generated $R$-module $M$.

Faltings [Fa1] established that the Annihilator Theorem for
cohomology modules holds over $R$ if $R$ is a homomorphic image of a
regular ring or $R$ has a dualizing complex. In [Ra2], Raghavan
proved that  the Local-global Principle for the annihilation of
local cohomology modules holds over $R$ if $R$ is a homomorphic
image of a regular ring. More recently, Khashyarmanesh and Salarian
[KS] obtained that the Annihilator Theorem and the Local-global
Principle for cohomology modules hold over a homomorphic image of a
(not necessarily finite-dimensional) Gorenstein ring.

Clearly, if the Annihilator Theorem holds over a ring $R$,  one can
use it to prove the existence of some annihilator $x_m$ of the local
cohomology modules $H_m^i(R)$ for each maximal ideal $m$. However,
the element $x_m$ may be dependent on the choice of the maximal
ideal m.     Raghavan [Ra1, Theorem 3.1] established a interesting
uniform annihilator theorem of local cohomology modules which states
that if $R$ is a homomorphic image of a biequidimensional regular
ring of finite dimension and $M$ a finitely generated $R$-module,
then there exists a positive integer $k$ (depending only on $M$)
such that for any ideals $I, J$ of $R$, we have
$J^k\text{H}_I^i(M)=0$  for $i< \lambda_I^J(M)$.

Note that, in all the results mentioned  above, the ring
considered must be, at least, a homomorphic image of a Gorenstein
ring. In this paper,  we first study the properties of the rings
containing uniform local cohomological annihilators. It turns out
that all these rings should be universally catenary and locally
equidimensional (Theorem 2.1). Due to this fact, we are able to
show that a power of a uniform local cohomological annihilator is
a strong uniform local cohomological annihilator (Theorem 2.2). We
will establish a useful criterion for the existence of uniform
local cohomological annihilators.  An easy consequence of one of
our main results shows that if a locally equidimensional
noetherian ring $R$ of positive dimension is a homomorphic image of a Cohen-Macaulay
(abbr. $CM$) ring of finite dimension or an excellent local ring,
then $R$ has a uniform local cohomological annihilator. This
greatly generalizes a lot of known results. Especially, it gives a
positive answer to a conjecture of Huneke [Hu, Conjecture 2.13] in
the local case.

 The  technique of the paper is   different from the technique
used by Roberts [Ro].  The point of  our technique is that  a
uniform local cohomological annihilator  of a ring $R$ may  be
chosen only dependent on the dimension of $R$ and the multiplicity
of each minimal prime ideal of $R$. One of our main results of the
paper is the following theorem,  which essentially reduces the
property that a ring $R$ has a  uniform local cohomological
annihilator to the same property for $R/P$ for all minimal prime
$P$ of $R$. Explicitly:

 \vspace* {0.3cm}

\noindent {\bf Theorem 3.2.}\ \ {\it Let $R$ be a
noetherian ring of finite dimension $d>0$. Then the following conditions
are equivalent:

{{\rm{(i)}}}\ \ $R$ has a   uniform local cohomological
annihilator.

{{\rm{(ii)}}}\  $R$ is locally equidimensional, and  $R/P$ has a
uniform local cohomological annihilator  for each minimal prime
ideal $P$ of $R$. }

 \vspace* {0.3cm}

In section 4, we discuss the uniform  local cohomological
annihilators for excellent rings. The main result of this section is
the following, which shows that the conjecture of Huneke [Hu,
Conjecture 2.13] is valid if the dimension of the ring considered is
no more than 5.

\vspace* {0.3cm}

\noindent {\bf Theorem 4.6.}\ \ {\it Let $R$ be a locally
equidimensional  excellent  ring  of \linebreak dimension
$d>0$. If $d\leq 5$, then $R$ has a uniform local cohomological
annihilator. }

\vspace* {0.3cm}

\section{ Basic properties}

 Now we begin with studying the properties of  a  noetherian
 ring $R$
containing a  uniform local cohomological annihilator. Quite
unexpectedly, it turns out that $R$ must be locally
equidimensional and universally catenary.

 Let $R$ be  a  noetherian ring. For a maximal ideal $m$ of
 $R$,   one can see easily from the definition of local cohomology that
$$\text{H}_{mR_m}^i(R_m) = \text{H}_m^i(R)$$
for $i\geq 0$. If $x$ is a  uniform local cohomological
annihilator of $R $, then $\text{H}^i_m(R)=0$ for every maximal
ideal $m$ with $x\notin m$ and $i< \text{ht}(m)$. So, for such a
maximal ideal $m$, $R_m$ is $CM$, and thus  $R_x$ is a $CM$ ring.

\begin{theorem} Let $R$ be a  noetherian and $x\in {R^\circ}$
 a  uniform local cohomological  annihilator  of $R $. Then

{\rm{(i)}} $R$ is locally equidimensional.

{\rm{(ii)}} $R$ is universally catenary.

\end{theorem}
\begin{proof} (i) \  Suppose that, on the contrary, $R$ is not locally
equidimensional. It implies  that there exists a maximal ideal $m$
of $R$, and a minimal prime ideal $P$  contained in $m$ such that
$\text{ht}(m/P) <\text{ht} (m) $.  Replacing $R$ by $R_m$, we can
assume $R$ is a local ring and $m$ is the unique maximal ideal of
$R$. Let
$$0 = q\cap q_2\cap \cdots \cap q_r$$
denote a shortest primary decomposition of the zero ideal of $R$,
where $q$ is $P$-primary.   By the choice of $P$ and $m$, it is
clear, $r>1$. So we can choose an element $y\notin P$ such that
$yq = 0$. Consequently $y\text{H}^i_m(q) = 0$ for $i\geq 0$.

  Consider the short exact sequence
$$0 \rightarrow  q \rightarrow R \rightarrow R/q\rightarrow 0.$$
It induces the following long exact sequence
$$\cdots \rightarrow  \text{H}_m^i(R) \rightarrow \text{H}_m^i(R/q)\rightarrow
\text{H}_m^{i+1}(q) \rightarrow \cdots.$$ Since $x\text{H}^i_m(R) =
0$ for $i<\text{ht} (m)$, we conclude $xy\text{H}^i_m(R/q) = 0$ for
$i<\text{ht} (m)$. In particular, we have $$xy\text{H}^e_m(R/q) = 0
\eqno (2.1)$$ where $e = \text{ht}(m/q)$. As $xy$ is a non
zero-divisor for $R/q$, the sequence
$$0\rightarrow  R/q\overset {xy}\to R/q\rightarrow
R/{(q+(xy))}\rightarrow 0$$ is a short exact sequence of
$R$-modules. So we have an exact sequence

$$\text{H}_m^e(R/q) \overset {xy}\to \text{H}_m^e(R/q)\rightarrow
\text{H}_m^{e}(R/(q+(xy))).$$
 Note that $\text{dim}(R/(q+(xy))) = e-1$. Hence  $\text{H}_m^{e}(R/(q+(xy))) = 0$ by [Gr, Proposition 6.4]. It shows the
morphism $\text{H}^e_m(R/q) \overset {xy}\to \text{H}^e_m(R/q)$ is
surjective. Thus $$\text{H}^e_m(R/q) = 0$$ by (2.1), but this is
impossible by [Gr, Proposition 6.4] again. Therefore $R$ is locally
equidimensional.

(ii) \ To prove the conclusion, it suffices to prove $R_m$ is
universally catenary for every maximal ideal $m$. So we can assume
that $R$ a local ring and $m$ is the unique maximal ideal of $R$.
By a theorem of Ratliff (see [Ma, Theorem 15.6]), it is enough to
prove $R/P$ is universally catenary for every minimal prime ideal
$P$.

For a fixed minimal prime ideal $P$, let
$$0 = q\cap q_2\cap \cdots \cap q_r$$
denote a shortest primary decomposition of the zero ideal of $R$,
where $q$ is $P$-primary. If $r = 1$, $x$ is a non zero-divisor of
$R$. Clearly, $x$ is also a  uniform local cohomological
annihilator  of $\hat{R}$, where $\hat{R}$ is the $m$-adic
completions of $R$. Hence by (i), $\hat{R}$ is equidimensional. It
follows from another theorem of Ratliff (see [Ma, Theorem 31.7])
that $R$ is universally catenary.

If $r>1$, choose an element $y\notin P$ as in the proof (i) such
that $xy$ is a non zero-divisor of $R/q$ and so the image of $xy$
in $R/q$ is a uniform local cohomological annihilator  of $R/q$.
Just as in the case $r = 1$, we assert that $R/q$ is universally
catenary, and consequently, $R/P$ is universally catenary. This
proves (ii).
\end{proof}

It is easy to see that a strong  uniform local cohomological
annihilator  of a noetherian $R$ is also a uniform local
cohomological annihilator. Conversely, we have:

\begin{theorem} Let $R$ be a  noetherian ring of finite  dimension $d$ and $x$
be a  uniform local cohomological  annihilator  of $R $. Then a
 power of $x$ is  a  strong uniform local  cohomological
annihilator  of $R $.
\end{theorem}

\begin{proof}

 First of all, we assume that $R$ is a unmixed local ring with the maximal
ideal $m$. Let $\hat{R}$ denote the completion of $R$ with respect
to $m$.  As $x$ is not contained in any minimal prime ideal of $R$
and $R$ is unmixed, $x$ must be a non zero-divisor of $R$, and so
$x$ is a non zero-divisor of $\hat{R}$. Hence the element $x$ is
not contained in any minimal prime ideal of $\hat{R}$. Since
$H_m^i(R)\simeq H_{m\hat{R}}^i(\hat{R})$ for all $i$, it shows
$xH_{m\hat{R}}^i(\hat{R})=0$ for $i<d$. Therefore we conclude that
$x$ is a uniform local cohomological annihilator of $\hat{R}$.

According to Cohen Structure Theorem for complete ring, one can
write $\hat{R}=S/I$, where $S$ is a Goreinstein local ring of
dimension d. By local duality, we have $x\text{Ext}_S^i(\hat{R},
S)=0$ for $i>0$. Hence
$$x\text{Ext}_{S_Q}^i(\hat{R}_Q, S_Q)=0$$
for every prime ideal $Q$ of $S$ and $i>0$. By local duality
again, we conclude that for every non minimal prime ideal
$\bar{Q}$ of $\hat{R}$,
$$x\text{H}_{\bar{Q}}^i(\hat{R}_{\bar{Q}})=0$$ for $i<
\text{ht}(\bar{Q})$.

Let $P$ be an arbitrary non minimal prime ideal of $R$. Set $
\text{ht}(P)=r$. We can choose elements $x_1, x_2, \cdots, x_r$
contained in $P$ such that $\text{ht}(x_1, x_2, \cdots, x_r)=r$. Set
$I=(x_1, x_2, \cdots, x_r)$. Since
dim$(\hat{R}/I\hat{R})=\text{dim}(R/I)$, it is clear that
dim$(\hat{R}/I\hat{R})=d-r$. By Theorem 2.1, $\hat{R}$ is
equidimensional, and so $\text{ht}(I\hat{R})=r$. It follows from
[HH1, Theorem 11.4 (b)] that for every $i\geq 1$ and for all $t>0$
$$x^{2^d-1}\text{H}_i(x_1^t, x_2^t, \cdots, x_r^t; \hat{R})=0$$
where $\text{H}_i(x_1^t, x_2^t, \cdots, x_r^t; \hat{R})$ denotes
the Koszul homology group.

Since $\hat{R}$ is faithfully flat over $R$, we have
$x^{2^d-1}\text{H}_i(x_1^t, x_2^t, \cdots, x_r^t; R)=0$ for every
$i\geq 1$ and for all $t>0$. Hence,  we obtain
$x^{2^d-1}\text{H}_I^i(R)=0$ for
 $i<r$. In particular
  $$x^{2^d-1}\text{H}_{PR_P}^i(R_P)=x^{2^d-1}(\text{H}_I^i(R))_P=0.$$
Therefore $x^{2^d-1}$ is a  strong uniform local cohomological
annihilator  of $R $.

Secondly, we assume that $R$ is  a  unmixed  ring. By the above
proof of the  local case, we have for any maximal ideal $m$ of $R$
and any prime ideal $P\subseteq m$,
$$x^{2^{{d_m}-1}}\text{H}^i_{PR_P}(R_P)=0$$ for $i<\text{ht}P$, where
$d_m$ stands for the dimension of the local ring $R_m$. Since
$d_m\leq d$ for every maximal ideal $m$, we conclude that
 $x^{2^d-1}$ is a strong uniform local
cohomological annihilator of $R $.

Finally, we assume that $R$ is not a  unmixed  ring.  Let $P_1, P_2,
\cdots, P_r$ be all the distinct minimal prime ideals of $R$. Let
$$0 = q_1\cap q_2\cap \cdots \cap q_r\cap q_{r+1}\cap \cdots \cap q_t$$
denote a shortest primary decomposition of the zero ideal of $R$,
where $q_i$ is $P_i$-primary for $i\leq r$. Clearly,  $t>r$. Set
$I=q_1\cap q_2\cap \cdots \cap q_r$, and $S=R/I$. Since $x$ is a
uniform local cohomological annihilator of $R$, it implies that
$\text{H}_m^i(R)=0$ for any maximal ideal $m$ with $x\notin m$ and
$i<\text{ht}(m)$, i.e. $R_x$ is $CM$. Hence $x$ must be contained in
every embedded associated primes of $R$. From this fact, one can
choose a positive integer $n_1$ such that $x^{n_1}\in q_j$ for all
$j$ with $r+1\leq j\leq t$. So   $x^{n_1}I=0$. Replacing $x$ by
$x^{n_1}$, we may assume that $xI=0$. Thus
$$x\text{H}_{PR_P}^i(I_P)=0 \eqno (2.2)$$ for all $i$ and for all prime
ideals $P$ of $R$.

Consider
 the short exact sequence
 $$0\rightarrow I\rightarrow R\rightarrow S\rightarrow 0. \eqno (2.3)$$
For any maximal ideal $m$ of $R$, we have the following long exact
sequence
$$\cdots\rightarrow \text {H}^{i}_{m}(R)\rightarrow
\text {H}^{i}_m(S)\rightarrow\text
{H}^{i+1}_{m}(I)\rightarrow\cdots. $$ By the definition of a uniform
local cohomological annihilator and (2.2), we conclude that
$x^2\text {H}^{i}_m(S)=0$ for $i<\text{ht}(m)$. Let $\bar {m}$ to
 denote the image of $m$ in $S$, we have
 $\text{ht}(m)=\text{ht}(\bar{m})$ by Theorem 2.1. So $x^2\text {H}^{i}_{\bar
{m}}(S)=0$ for  any maximal ideal $\bar{m}$ of $S$ and for $i<
\text{ht}(\bar{m})$. Hence the image
 of $x^2$ in $S$ is a uniform local cohomological annihilator of
 $S$. Clearly, $S$ is unmixed.  By the unmixed case proved above,
 we assert that there exists a positive integer $n_2$ such that
the image of  $x^{2n_2}$ in $S$ is a strong uniform local
cohomological annihilator of
 $S$. Replacing $x$
by  $x^{2n_2}$, we may assume that the image of $x$ in $S$ is a
strong uniform local cohomological annihilator of
 $S$.

 Let $P$ be an arbitrary prime ideal of $R$. We use $\bar {P}$ to
 denote the image of $P$ in $S$. By Theorem 2.1, we have
 $\text{ht}(P)=\text{ht}(\bar{P})$. Hence by the choice of $x$, we
 conclude that $$xH_{PR_P}^i(S_P)=0 \eqno(2.4)$$ for $i< \text{ht}(P)$.
 Localizing the short exact sequence (2.3) at $P$, it induces
 the following long exact sequence
$$\cdots\rightarrow \text {H}^{i}_{PR_P}(I_P)\rightarrow
\text {H}^{i}_{PR_P}(R_P)\rightarrow\text
{H}^{i}_{PR_P}(S_P)\rightarrow\cdots. $$ Hence it follows from (2.2)
and (2.4) that $x^2\text {H}^{i}_{PR_P}(R_P)=0$ for $i<
\text{ht}(P)$. Therefore $x^2$ is a strong uniform local
cohomological annihilator of $R$, and this ends the proof of the
theorem.

\end{proof}

\section{ Equivalent conditions }

In this section, we will prove one of our main result, which
essentially reduces the property that a ring $R$ has a uniform
local cohomological annihilator to proving the same property for
$R/P$ for all minimal primes of $R$. This reducing process is very
useful, it enable us to find a  uniform local cohomological
annihilator more easily and directly. Now, before we prove the
main result of this section, we need a lemma which will play a key
role in the rest of the section.

\begin{lemma} Let $(R, m)$ be a  noetherian local ring
 of dimension $d$ , and $P$ be a minimal prime ideal of $R$.
 Let $$0\rightarrow R/P\rightarrow R\rightarrow N_1\rightarrow 0$$
$$0\rightarrow R/P\rightarrow N_1\rightarrow N_2\rightarrow 0     \eqno (3.1)$$
$$       \cdots        $$
$$0\rightarrow R/P\rightarrow N_{t-1}\rightarrow N_t\rightarrow 0$$
be a series of short exact  sequences of finitely generated
$R$-modules. Let $y $ be an element of $R$ such that $yN_t = 0$.

{\rm{(i)}}  \ \ If there is an element $x$ of $R$ such that
$x\text{H}^i_m(R) = 0$ for $i<d$, then
$(xy)^{t^{d-1}}\text{H}^i_m(R/P) = 0$   for $i<d$.

{\rm{(ii)}}  \ \ If there is an element $x$ of $R$ such that
$x\text{H}^i_m(R/P) = 0$ for $i<d$, then $x^ty\text{H}^i_m(R) = 0$
for $i<d$.

\end{lemma}

\begin{proof} Since $R$ is , we have $d>0$.

(i)\ \ By the choice of $y$, it implies
 $$y\text{H}^i_m(N_t) = 0      \eqno (3.2).$$
for all $i\geq 0$. We will use induction on $i$ to prove
$$(xy)^{t^{i}}\text {H}^i_m(R/P) = 0 $$
holds for  $ 0\leq i< d$.

For $i = 0$, it is trivial because $\text {H}^0_m(R/P) = 0$. Now,
for $0< i< d$, Suppose that we have proved
$$(xy)^{t^{i-1}}\text {H}^{i-1}_m(R/P) = 0.  \eqno (3.3)$$

Set $k = t^{i-1}$. Let us consider the long exact sequence of
local cohomology derived from the last short exact sequence in
(3.1):

$$\cdots\rightarrow \text {H}^{i-1}_m(R/P)\rightarrow
\text {H}^{i-1}_m(N_{t-1})\rightarrow\text
{H}^{i-1}_m(N_t)\rightarrow\cdots. $$   By (3.2) and (3.3), it
follows ${(xy)}^ky\text {H}^{i-1}_m(N_{t-1}) = 0$. Continue the
process, one can prove $(xy)^{jk}y\text {H}^{i-1}_m(N_{t-j}) = 0$
for $j = 1, 2, \cdots, t-1$. Hence by the long exact sequence of
local cohomology derived from the first short exact sequence in
(3.1):
$$ \cdots\rightarrow \text {H}^{i-1}_m(N_1)\rightarrow
\text {H}^i_m(R/P)\rightarrow\text {H}^i_m(R)\rightarrow\cdots $$
and the condition $x\text{H}^i_m(R) = 0$, we have
$(xy)^{(t-1)k+1}\text {H}^i_m(R/P) = 0$. It easy to check,
$t^{i}\geq (t-1)k+1$, so it follows that $(xy)^{t^{i}}\text
{H}^i_m(R/P) = 0$. This completes the inductive proof. In
particular, we have proved
 $(xy)^{t^{d-1}}\text{H}^i_m(R/P) = 0$ for $i<
d$.

(ii) By the condition, we have
 $$x\text{H}^i_m(R/P) = 0  \eqno (3.4)$$
 for $i<d$. Set $R = N_0$.

   Now, we will use induction on $j$ to  prove that
$$x^jy\text{H}^i_m(N_{t-j}) = 0 \ \ \text{for}\ \
i<d  $$ hold for $0\leq j\leq t$.

For $j = 0$, it is trivial by (3.2). For $j>0$, suppose that we
have proved
$$x^{j-1}y\text{H}^i_m(N_{t-(j-1)}) = 0  \eqno (3.5)$$ for
$i<d$.

 Consider the $t-(j-1)$-th short exact
sequence in (3.1), it induces the following long exact sequence
$$\cdots\rightarrow \text {H}^i_m(R/P)\rightarrow
\text {H}^i_m(N_{t-j})\rightarrow\text
{H}^i_m(N_{t-(j-1)})\rightarrow\cdots.$$
 By (3.4) and (3.5), we conclude
$$x^jy\text{H}^i_m(N_{t-j}) = 0$$ for
$i<d$. This completes the inductive proof. In particular, we have
$x^ty\text{H}^i_m(R) = 0$ for $i<d$.

\end{proof}

We are now ready to prove our main result of the section.
\begin{theorem} Let $R$ be a  noetherian ring of finite  dimension $d>0$.
Then the following conditions are equivalent:

{{\rm{(i)}}} \  $R$ has a  uniform local cohomological
annihilator.

{{\rm{(ii)}}}\   $R$ is locally equidimensional, and $R/P$ has a
uniform local cohomological annihilator  for each minimal prime
ideal $P$ of $R$.

\end{theorem}

\begin{proof}

 (i) $\Rightarrow$ (ii) \ \  The first conclusion of (ii) comes from Theorem 2.1.
 Let $P$ be an arbitrary minimal prime
ideal of $R$.  Let $x$ be a uniform local cohomological
annihilator  of $R$. Put $t = l(R_P)$. It is easy to see that
there exist finitely generated $R$-modules $N_1, N_2, \cdots, N_t$
such that

(1)\ \ $N_1, N_2, \cdots, N_t$ fit into a series of the following
short exact sequences

$$0\rightarrow R/P\rightarrow R\rightarrow N_1\rightarrow 0$$
$$0\rightarrow R/P\rightarrow N_1\rightarrow N_2\rightarrow 0     \eqno (3.6)$$
$$       \cdots        $$
$$0\rightarrow R/P\rightarrow N_{t-1}\rightarrow N_t\rightarrow 0$$

(2)\ there exists an element $y\in R\setminus P$, $yN_t = 0.$

Clearly, $(xy)^{t^{d-1}}$ is not contained in $P$. We will show
that the image of $(xy)^{t^{d-1}}$ in $R/P$ is a  uniform local
cohomological annihilator  of $R/P$. Since $R$  is locally
equidimensional,, we have $\text{ht}(m/P)=\text{ht}(m)$ for every
maximal ideal $m$ with $m\supseteq P$. Thus it suffices to prove
that, for every maximal ideal $m$ with $m\supseteq P$

$$(xy)^{t^{d-1}}\text{H}_{mR_m}^i((R/P)_m)  =  0$$
for $i<\text{ht}(m)$.

Localizing the short exact sequences in (3.6) at $m$, we obtain
the following the short exact sequences

$$0\rightarrow (R/P)_m\rightarrow R_m\rightarrow (N_1)_m\rightarrow 0$$
$$0\rightarrow (R/P)_m\rightarrow (N_1)_m\rightarrow (N_2)_m\rightarrow 0     $$
$$       \cdots        $$
$$0\rightarrow (R/P)_m\rightarrow (N_{t-1})_m\rightarrow (N_t)_m\rightarrow 0.$$

By the choice of $x$, we have $x\text{H}_{mR_m}^i(R_m)  = 0$ for
$i<\text{ht}(m)$. Clearly, $y(N_t)_m = 0$. Hence by Lemma 3.1 (i),
we conclude $(xy)^{t^{e-1}}\text{H}_{mR_m}^i((R/P)_m)  =  0$ for
$i<e$, where $e= \text{ht}(m)$. Therefore for every maximal ideal
$m$ with $m\supseteq P$, we have
$$(xy)^{t^{d-1}}\text{H}_{mR_m}^i((R/P)_m)  =  0$$ for
$i<\text{ht}(m)$. Hence the image of $(xy)^{t^{d-1}}$ in $R/P$ is
a  uniform local cohomological  annihilator  of $R/P$,   and this
proves  (i) $\Rightarrow$ (ii).

(ii) $\Rightarrow$ (i) \ \ Let $P_1, P_2, \cdots, P_r$ be all the
distinct minimal prime ideals of $R$.  For each $j$ $1\leq j\leq r$,
Put $l(R_{P_j})  =  t_j$.  For a fixed $j$, It is easy to see that
there exist finitely generated $R$-modules $N^{(j)}_1, N{(j)}_2,
\cdots, N^{(j)}_{t_j}$  satisfying the following two properties:

(1)\ \ $N^{(j)}_1, N^{(j)}_2, \cdots, N^{(j)}_{t_j}$ fit into a
series of the following short exact sequences

$$0\rightarrow R/P_j\rightarrow R\rightarrow N^{(j)}_1\rightarrow 0$$
$$0\rightarrow R/P_j\rightarrow N^{(j)}_1\rightarrow N^{(j)}_2\rightarrow 0     \eqno (3.7)$$
$$       \cdots        $$
$$0\rightarrow R/P_j\rightarrow N^{(j)}_{t_j-1}\rightarrow N^{(j)}_{t_j}\rightarrow 0$$

(2)\ \ There exists an element $y_j\in R\setminus P$ such that
 $y_jN_{t_j} = 0.$ and $y_j$ lies in all $P_k$ except $P_j$.

To prove the conclusion, it is enough to find an element $x\in
R^\circ$, such that for every maximal ideal $m$
$$x\text{H}^i_{mR_m}(R_m) = 0$$
for $i< \text{ht}(m)$.

By the condition,  for each $j$, there exists an element
$x_j\notin P_j$ such that its image in $R/P_j$ is a uniform local
cohomological annihilator  of $R/P_j$. Set $x = \sum
x_j^{t_j}y_j$. It is easy to check $x$ lies in no minimal prime
ideal of $R$. We will prove that $x$ is a  uniform local
cohomological annihilator  of $R$.

Let $m$ be an arbitrary non minimal prime ideal of $R$. Put $e =
\text{ht}(m)$. For a fixed $j$, localizing the short exact
sequences (3.7) at $m$, we obtain the following short exact
sequences
$$0\rightarrow (R/P_j)_m\rightarrow R_m\rightarrow (N^{(j)}_1)_m\rightarrow 0$$
$$0\rightarrow (R/P_j)_m\rightarrow (N^{(j)}_1)_m\rightarrow (N^{(j)}_2)_m\rightarrow 0     $$
$$       \cdots        $$
$$0\rightarrow (R/P_j)_m\rightarrow( N^{(j)}_{t_j-1})_m\rightarrow (N^{(j)}_{t_j})_m\rightarrow 0$$
and it is clear $y_j(N^{(j)}_{t_j})_m = 0$.

  If $m$ contains $P_j$,  we have
 $\text{ht}(m/P_j)=e$ by the assumption that $R$ is locally equidimensional.  Thus by the choice of $x_j$,
$$x_j\text{H}^i_{mR_m}((R/P_j)_m) = 0   \eqno (3.8) $$
for $i<e $. If $m$ does not contain $P_j$, the statement (3.8)
holds trivially.
 Hence by Lemma 2.2 (ii),  we conclude that
$$x^{t_j}y_j\text{H}^i_{mR_m}(R_m) = 0 \eqno (3.9) $$
for $i< e$.  Therefore, by the choice of $x$, we have

$$x\text{H}^i_{mR_m}(R_m) = 0, \ \ \text{for}\ \  i< \text{ht}(m) $$
holds for every maximal ideal $m$ of $R$. So $x$ is a uniform
local cohomological  annihilator  of $R$. This proves (ii)
$\Rightarrow$ (i).
\end{proof}

Now, we end this section by   an  interesting corollary, which is
not known in such a extent even in the local case.

\begin{corollary} Let $R$ be a locally
equidimensional  noetherian ring of finite positive dimension. Then $R$ has  a uniform
local cohomological annihilator in any one of the following cases:

{\rm{(i)}}  $R$ is the homomorphic image of a $CM$ ring of finite
 dimension.

{\rm{(ii)}} $R$ is an excellent local ring.
\end{corollary}
\begin{proof} (i)\ Represent $R$ as $R = S/I$, where $S$ is a
$CM$ ring of finite  dimension and $I$ an ideal of $S$. By Theorem
3.2, it suffices to prove that for any minimal prime ideal $P$ of
$R$, $R/P$ has  a  uniform local cohomological annihilator. Let
$Q$ be the prime ideal of $S$ such that $P = Q/I$. It suffices to
prove that $S/Q$ has  a  uniform local cohomological annihilator.

Set $n = \text{ht}(Q)$. If $n = 0$, then $Q$ is a minimal prime
ideal of $S$. Note that $1$  is  a uniform local cohomological
annihilator  of $S$,   the conclusion follows immediately from
Theorem 3.2.

Assume that $n>0$. Since $S$ is $CM$, we can choose a regular
sequence $x_1, x_2, \cdots , x_n$ contained in Q. Thus $S/(x_1,
x_2, \cdots , x_n)$ is still a $CM$ ring. It is obvious that
$Q/(x_1, x_2, \cdots , x_n)$ is a minimal prime ideal of $S/(x_1,
x_2, \cdots , x_n)$. Thus from the case $n = 0$, we assert $S/Q$
has  a uniform local cohomological annihilator by Theorem 3.2.

(ii)\  By [HH2, Lemma 3.2], every excellent local domain has  a
strong uniform local cohomological annihilator,  and thus $R$ has
a strong uniform local cohomological annihilator by Theorem 3.2.
\end{proof}

\section{Uniform local cohomological annihilators of excellent rings   }

In this section, we restrict our discussion to the uniform local
cohomological annihilators of excellent rings. By two theorems of
Hochster and   Huneke [HH1, Theorem 11.3, 11.4], it is easy to see
that
  Huneke's conjecture [Hu, Conjecture 2.13] is valid if every
locally equidimensional excellent noetherian ring $R$ of finite
dimension has a strong uniform local cohomological annihilator.
Since for each positive integer $n$, the Koszul complex of every
sequence $x_1^n, x_2^n, \cdots, x_k^n$ in $R$ with ht$(x_1,
x_2,\cdots, x_k)=k$ is a complex satisfies the standard conditions
on  height and rank in the sense of [Hu, (2.11)],  it is easy to
see that the converse of this result is also true. Due to Theorem
2.2 and Theorem 3.2,  Huneke's conjecture is equivalent to the
following:

\begin{conjecture}
Let $R$ be an  excellent noetherian domain of finite dimension.
Then $R$ has a uniform local cohomological annihilator.
\end{conjecture}
 The conjecture is known to be true if $R$ is an excellent
 normal domain  of dimension $d\leq 3$ [Hu,Proposition 4.5(vii)].
 We will prove that Conjecture 4.1 is true for the
ring $R$ with $\text{dim}(R)\leq 5$.  In order to prove the main
result of this section, we need the following  result established
by Goto[Go, Theorem 1.1].

\begin{proposition}
Let $(R, m)$ be a local ring of dimension $d$ and $x$ is an
element of $m$ with $(0 : x)=(0 : x^2)$. Then the following
conditions are equivalent.

{\rm{(i)}} $R/{x^nR}$ is a $CM$ ring of dimension $d-1$ for every
integer $n>0$.

{\rm{(ii)}} $R/{x^2R}$ is a $CM$ ring of dimension $d-1$.

\end{proposition}

Another important result we need is the following explicit version
of \linebreak Corollary 3.3 (ii), which follows  from [HH2, Lemma
3.2].

\begin{proposition}
Let $(R, m)$ be an excellent  local domain of dimension $d>0$ and
$x$ be an element of $m$ such that $R_x$ is a $CM$ ring. Then a
power of $x$ is a uniform local cohomological annihilator of $R$.

\end{proposition}

\begin{proof}

Let $x_1, x_2,\cdots, x_d$ be an arbitrary system of parameters in
$m$. \linebreak By    [HH2, Lemma 3.2], there exists a positive
integer $n$ such that $x^n$ kills all the higher Koszul homology
$\text{H}_i(x_1^{n_1}, x_2^{n_2}, \cdots, x_d^{n_d}, R), i>0$ for
all positive integers $n_1, n_2, \cdots n_d$. Hence

$$x^n\text{H}_m^i(R)=\underset{t\rightarrow\infty}{\lim }\text{H}_{d-i}(x_1^{t},
x_2^{t}, \cdots, x_d^{t}, R)=0  $$ for $i<d$. It shows that $x^n$
is a uniform local cohomological annihilator of $R$.
\end{proof}

Let $R$ be an excellent ring of dimension $d>0$. For any prime
ideal $P$ of $R$, the regular locus of $R/P$ is a non-empty open
subset of Spec$(R/P)$, and so there exists a non-empty open subset
$U$ of Spec$(R/P)$ such that for any $Q\in U$, $(R/P)_Q$ is $CM$.
By Nagata Criterion for openness, we conclude that the $CM$ locus
of $R$ is  open in Spec$(R)$ (see [Ma, Theorem 24.5]). Moreover,
as every minimal prime ideal $P$ of $R$ lies in the $CM$ locus of
$R$, we assert that the $CM$ locus of $R$ is a non trivial open
set in Spec$(R)$. Hence we can choose an element $x\in R^\circ$
such that $R_x$ is a $CM$ ring. Hence we
can choose an element $x\in R^\circ$ such that $R_x$ is a $CM$ ring.
By Proposition 4.3, for any maximal ideal $m$ of $R$, there exists a
positive integer $n_m$ such that $x^{n_m}\text{H}_m^i(R)=0$ for
$i<\text{ht}(m)$. Clearly, the positive integer $n_m$ may be
dependent on $m$. To solve Conjecture 4.1, it suffices to find a
positive integer $n$ such that it is independent on the choices of
the maximal ideals $m$. Nevertheless, we have the following useful
corollary.

\begin{corollary}
Let $R$ be an excellent   domain of dimension $d>0$ and $x$ be an
element of $m$ such that $R_x$ is a $CM$ ring. If $T$ is a finite
set of maximal ideals of $R$, then there exists a positive integer
$n$ such that for any $m\in T$, $x^n\text{H}_m^i(R)=0$ for
$i<\text{ht}(m)$.

\end{corollary}

To simplify the proof of  the main result of this section, we need
the following lemma, which enable us to obtain the annihilators of
local cohomology modules.

\begin{lemma}
Let $(R, m)$ be a noetherian local ring of dimension $d$. Let
$x_1, x_2, \cdots, x_r$ be a part of system of parameters in $m$
and $x$ an element in $m$. Suppose that

(\rm i)\ \ $R/(x_1^{n_1}, x_2^{n_2}, \cdots, x_r^{n_r})$ are $CM$;

(\rm{ii})\ \ For $1\leq i\leq r$, $x((x_1^{n_1}, x_2^{n_2},
\cdots, x_{i-1}^{n_{i-1}}) : x_i^{n_i})\subset (x_1^{n_1},
x_2^{n_2}, \cdots, x_{i-1}^{n_{i-1}})$

\noindent hold for all positive integers $n_1, n_2, \cdots, n_r$.
Then $x^{r}H^i_m(R)=0$ for $i < d$.

\end{lemma}

\begin{proof}

we will use induction on $j (0\leq j\leq r)$ to assert that for
any positive integers $n_1, n_2, \cdots, n_j$
$$x^{r-j}(\text{H}_m^i(R/(x_1^{n_1}, x_2^{n_2}, \cdots, x_j^{n_j}))=0$$
for $i< d-j$, and then the lemma follows if we set $j=0$.

By the assumption, for  arbitrary fixed integers $n_1, n_2,
\cdots, n_r$, \linebreak $R/(x_1^{n_1}, x_2^{n_2}, \cdots,
x_r^{n_r})$ is $CM$, so the conclusion is trivial in this case.
Suppose that we have proved the conclusion for $t+1\leq j\leq r$.
For a fixed $i$ with $i<d-t$, let $z$ be an arbitrary element in
$\text{H}_m^i(R/(x_1^{n_1}, x_2^{n_2}, \cdots, x_t^{n_t}))$.
Choose a positive integer $n_{t+1}$ such that
$x_{t+1}^{n_{t+1}}z=0$. Put

(1)  $R_t=R/(x_1^{n_1}, x_2^{n_2}, \cdots, x_t^{n_t})$

(2)  $R_{t+1}=R/(x_1^{n_1}, x_2^{n_2}, \cdots, x_{t+1}^{n_{t+1}})$

(3) $U_{t}=((x_1^{n_1}, x_2^{n_2}, \cdots, x_{t}^{n_{t}}) :
x_{t+1}^{n_{t+1}})/(x_1^{n_1}, x_2^{n_2}, \cdots, x_{t}^{n_{t}})$ ;

(4) $N_{t}=x_{t+1}^{n_{t+1}}(R/(x_1^{n_1}, x_2^{n_2}, \cdots,
x_t^{n_{t}}))$ .

Let us consider the  short exact sequences
$$0\rightarrow N_t\rightarrow R_t\rightarrow
R_{t+1}\rightarrow 0$$

$$0\rightarrow U_{t}\rightarrow R_{t}\rightarrow N_{t} \rightarrow 0$$
we have the following long exact sequences of local cohomology
$$\cdots\rightarrow \text{H}^{i-1}_m(R_{t+1})\rightarrow
\text{H}^{i}_m(N_{t})\overset {\phi_i}\to
\text{H}^{i}_m(R_{t})\rightarrow \cdots$$

$$\cdots\rightarrow \text{H}^{i}_m(U_{t})\rightarrow
\text{H}^{i}_m(R_{t})\overset {\psi_i}\to
\text{H}^{i}_m(N_{t})\rightarrow \cdots.$$ It is easy to see that
the composition $\phi_i \psi_i$ of $\psi_i$ and $\phi_i$ is the
morphism
$$\text{H}^{i}_m(R_{t})\overset{x_{t+1}^{n_{t+1}}}\longrightarrow
\text{H}^{i}_m(R_{t})$$ and thus $\phi_i (\psi_i(z))=0$. By the
induction hypothesis, $ x^{r-t-1}(\text{H}_m^{i-1}(R_{t+1}))=0$,
so from the first long exact sequence above, we conclude that
$(\psi_i(x^{r-t-1}z))=0$. The condition (ii) implies that
$xU_t=0$,  so $x\text{H}^{i}_m(U_{t})=0$. Hence from the second
long exact sequence above, we have $x^{r-t}z=0$. By the choice of
$z$, we have proved $x^{r-t}(\text{H}_m^{i}(R_{t}))=0$. This ends
the inductive proof of the lemma.

\end{proof}

In the rest of the paper, we will make use of Proposition 4.2,
Corollary 4.4 and Lemma 4.5 to prove the following main result of
this section.
\begin{theorem} Let $R$ be a locally equidimensional  excellent  ring  of dimension $d>0$. If $d\leq 5$, then $R$ has
a uniform local cohomological annihilator.
\end{theorem}

\begin{proof} By Theorem 3.2, we  may  assume that $R$ is a
excellent domain. Moreover, let $S$ be the integral closure of $R$ in the
field of fractions of $R$. Since $R$ is excellent, it follows that
$S$ is a finitely generated $R$-module. We first conclude that if
$S$ has a uniform local cohomological annihilator $y$, then $R$
also has a uniform local cohomological annihilator.

In fact, as $y$ is integral over $R$, we have
$$y^s + a_1y^{s-1}+ \cdots +a_s=0$$
for some elements  $a_1,  \cdots, a_s$ contained in $R$ with
$a_s\neq 0$. So it is clear $a_s$ is a uniform local cohomological
annihilator  of $S$. Consider the following natural short exact
sequence

$$0\rightarrow R\rightarrow S\rightarrow M\rightarrow 0. \eqno (4.1)$$
where $M$ is  a finitely generated $R$-module.  It is easy to find
a nonzero element $x$ of $R$ such that $xM=0$. Now, for an
arbitrary maximal $m$ of $R$, the all the minimal prime ideals
$Q_1, Q_2, \cdots, Q_t$ of $mS$ are maximal. Thus by the
Mayer-Vietoris sequence  of local cohomology, we conclude that

$$\text {H}_m^i(S)\simeq \text {H}_{Q_1}^i(S)\oplus \cdots \oplus \text
{H}_{Q_t}^i(S).$$ and consequently $a_s\text {H}_m^i(S)=0$ for
$i<\text{ht}(m)$. From the long exact sequence of local cohomology
induced from (4.1), we have $a_sx\text {H}_m^i(R)=0$ for
$i<\text{ht}(m)$. Thus $a_s x$ is a uniform local cohomological
annihilator  of  $R$. So in the following proof we assume $R$ is a
normal domain.

If $d=2$, then $R$ is $CM$, and there is nothing to prove. So we
assume $d>2$. Set $V=\{P\in \text{Spec}(R)\mid R_P \ \ \text {is
not a} \ \ CM \ \ \text{ring}\}$. As $R$ is a excellent ring, the
$CM$ locus of $R$ is open in $\text{Spec}(R)$, so there exists an
ideal $I$ of $R$ such that $V=V(I)$. Clearly, $\text{ht}(I)\geq
3$. Choose elements $x_1, x_2, x_3$ contained in $I$ such that
$\text{ht}((x_i, x_j))=2$ for $i\neq j$, and $\text{ht}((x_1, x_2,
x_3))=3$.

\begin{claim}{1.}
The union $T_1$ of the sets of  associated prime ideals
$\text{Ass}_R(R/{(x_1^{n_1}, x_2^{n_2})})$ is a finite set, where
the union is taken over all positive integers $n_1, n_2.$
\end{claim}
\begin{proof}
 It is well known the union $T_0$ of the sets of  associated prime ideals
$\text{Ass}_R(R/{(x_1, x_2^{n_2})})$ is a finite set, where the
union is taken over all positive integers $n_2$. We assert that
$T_0=T_1$. By induction on $n_1$, this follows easily  from the
short exact sequences

$$0\rightarrow R/(x_1, x_2^{n_2})\rightarrow R/{(x_1^{n_1}, x_2^{n_2})}\rightarrow R/{(x_1^{n_1-1},
x_2^{n_2})} \rightarrow 0$$
 for $n_1>1$. This proves the claim.
\end{proof}

Let $T_2$  denote the union of $T_1$ and  the set of minimal prime
ideals of $(x_1, x_2, x_3)$. Clearly, $T_2$ is a finite set. Since
$ R_{x_3} $ is $CM$, by [HH2, Lemma 3.2], we can choose a positive
number $n$ such that, for all positive integers $n_1, n_2$, the
following hold
$$((x_1^{n_1}, x_2^{n_2}): x_3^{n})_P= ((x_1^{n_1}, x_2^{n_2}):
x_3^{n+1})_P   \eqno(4.2)$$ for every prime ideal $P$ lies in the
set $T_2$.

For any prime ideal $P$ with $(x_1, x_2, x_3)\subseteq P$ and
$P\notin T_2$, it is clear $\text{ht}(P)\geq 4$ and $\text {depth}
(R_P/{(x_1^{n_1}, x_2^{n_2}})_P)\geq 1$. Localizing at $P$ if
necessary, we may assume $P$ is a maximal ideal. Giving an element
$c\in R$ satisfying $Pc\subset ((x_1^{n_1}, x_2^{n_2}): x_3^{n})$,
it follows $Px_3^nc\subset (x_1^{n_1}, x_2^{n_2})$. By the choice
of $P$, we have $x_3^nc\in (x_1^{n_1}, x_2^{n_2})$, and so we
conclude that $P$ is not an associated prime of $R/{((x_1^{n_1},
x_2^{n_2}): x_3^n)}$. Therefore the associated prime ideals
$R/{((x_1^{n_1}, x_2^{n_2}): x_3^n)}$ are contained in $T_2$. Thus
by (4.2), we conclude that
$$((x_1^{n_1}, x_2^{n_2}): x_3^{n})= ((x_1^{n_1}, x_2^{n_2}):
x_3^{n+1}) $$ hold for all positive integers  $n_1, n_2$.

Similarly,   enlarging $n$  if necessary,  one can prove
$$((x_2^{n_2}, x_3^{n_3}): x_1^{n})= ((x_2^{n_2}, x_3^{n_3}):
x_1^{n+1}) $$

$$((x_1^{n_1}, x_3^{n_3}): x_2^{n})= ((x_1^{n_1}, x_3^{n_3}):
x_2^{n+1}). $$ Replacing $x_1, x_2, x_3$ by $x_1^{n}, x_2^{n},
x_3^{n}$ respectively, we have proved

$$((x_1^{n_1}, x_2^{n_2}): x_3)= ((x_1^{n_1}, x_2^{n_2}):
x_3^{2}) $$

$$((x_2^{n_2}, x_3^{n_3}): x_1)= ((x_2^{n_2}, x_3^{n_3}):
x_1^{2})  \eqno(4.3)$$

$$((x_1^{n_1}, x_3^{n_3}): x_2)= ((x_1^{n_1}, x_3^{n_3}):
x_2^2) $$ hold for all positive integers  $n_1, n_2, n_3$.

If  $\text{dim}(R/{(x_1, x_2, x_3)})=0$, then here are only finite
number of maximal ideals $m$ such that $R_m$ may not be a $CM$
ring. So by Corollary 4.4, a power of $x_1$ is a uniform local
cohomological annihilator of $R$. In particular, the conclusion of
the theorem holds for $d=3$.

Now, in the following we assume $\text{dim}(R/{(x_1, x_2,
x_3)})>0$. Since $R$ is excellent, we can choose an element $x_4$
such that $x_4$ is not contained in any minimal prime ideal of
$(x_1, x_2, x_3)$ and $(R/{(x_1^2, x_2^2, x_3^2)})_P$ is $CM$ for
all prime ideals $P$ with $x_4\notin P$. For such a prime ideal
$P$, we conclude by  Proposition 4.2 that $(R/{(x_1^{2}, x_2^{2},
x_3^{n_3})})_P$ are all $CM$ rings for all positive integers
$n_3$. By (4.3) and by  Proposition 4.2 again, we conclude that
$(R/{(x_1^{n_1}, x_2^{n_2}, x_3^{n_3})})_P$ are all $CM$ rings for
all positive integers $n_1, n_2, n_3$.

\begin{claim}{2.}
Let $m$ be a maximal ideal of $R $ such that $m$ does not contain
the ideal $(x_1, x_2, x_3, x_4)$. Then for every $x_j$ ($1\leq
j\leq 3$), $x_j^3\text{H}_m^i(R)=0$ for $i< \text{ht}(m)$.
\end{claim}
\begin{proof}
Let $m$ be an arbitrary maximal ideal with $x_4\notin m$. If one
of $x_1, x_2, x_3$ does not lie in $m$, then $R_m$ is $CM$, and
the conclusion is trivial if one of them does not lie in $m$. So
we may assume $m\supseteq (x_1, x_2, x_3)$. Clearly, we may also
assume that $R$ is a $d$-dimensional  local ring with the unique
maximal ideal $m$.

 Note that $R/(x_1^{n_1}, x_2^{n_2}, x_3^{n_3})$ is $CM$ of dimension $d-3$
 for any positive integers $n_1, n_2, n_3$. Moreover, $x_1, x_2$
 is a regular sequence in $R$, so together with (4.3), we have for
$1\leq i\leq 3$, $x_j((x_1^{n_1},  \cdots, x_{i-1}^{n_{i-1}}) :
x_i^{n_i})\subset (x_1^{n_1}, \cdots, x_{i-1}^{n_{i-1}})$ for each
j with $1\leq j\leq 3$. Hence the conclusion of the claim follows
immediately from Lemma 4.5, and this proves the claim.

\end{proof}

By Claim 2,  we have proved that if a maximal ideal $m$ does not
contain the ideal $(x_1, x_2, x_3, x_4)$, then for each $j$
($1\leq j\leq 3$), $x_1^3H_m^i(R)=0$ for $i<\text{ht}(m)$. Now,
Suppose that $\text{dim}(R/{(x_1, x_2, x_3, x_4)})=0.$  There are
only a finite number of maximal ideals $m$ with  $(x_1, x_2, x_3,
x_4)\subseteq m$.  So by Claim 2 and Corollary 4.4, a power of
$x_j$ ($1\leq j\leq 3$) is a uniform local cohomological
annihilator  of $R$.  In particular, the conclusion of the theorem
holds for $d=4$. In the following proof, we assume $\text
{dim}(R/{(x_1, x_2, x_3, x_4)})>0$.

\begin{claim}{3.}
Let $T_3$ be the set of all associated prime ideals of
$R/{(x_1^{i_1}, x_2^{i_2},x_3^{i_3})}$ for all positive integers
$i_1, i_2, i_3 $ satisfying $i_1+ i_2+ i_3\leq 6 $. Then for any
positive integers $n_1, n_2, n_3$, every associated prime ideal of
$R/{(x_1^{n_1}, x_2^{n_2},x_3^{n_3})}$ lies in $T_3$.

\end{claim}
\begin{proof}
We prove the conclusion by  induction on $n=n_1+ n_2+ n_3$.
Clearly, if $n\leq 6$, the conclusion holds by the assumption. In
the following we assume  $n> 6$. Suppose
 we have proved the conclusion  for $n_1+ n_2+ n_3<n$. It is easy to see there  exists
one of $n_1, n_2, n_3 $, say $n_3$, such that $n_3\geq 3$.

Suppose that $P$ is an associated prime ideal of $R/{(x_1^{n_1},
x_2^{n_2},x_3^{n_3})}$ and $P\notin T_3$. Localizing at $P$ if
necessary, we may assume that $R$ is a local ring with the unique
maximal ideal $P$. Let $c\notin (x_1^{n_1}, x_2^{n_2},x_3^{n_3})$
be an element satisfying

$$Pc\subset (x_1^{n_1}, x_2^{n_2},x_3^{n_3}).$$
Then we can express $c=x_1^{n_1}c_1+ x_2^{n_2}c_2+x_3^{n_3-1}c_3$
by the induction hypothesis. For an arbitrary element $z\in P$, it
implies that there exists $c_4\in R$ such that
$$x_3^{n_3-1}(zc_3-x_3c_4 )\in (x_1^{n_1},
x_2^{n_2}).$$ By (4.3), we have $zx_3c_3\in (x_1^{n_1}, x_2^{n_2},
x_3^2)$. Hence $Px_3c_3\subset (x_1^{n_1}, x_2^{n_2}, x_3^2)$.
Since $n_1+n_2+2<n$,  we conclude $x_3c_3\in (x_1^{n_1},
x_2^{n_2}, x_3^2)$ by the induction hypothesis again,.
Consequently $c\in (x_1^{n_1}, x_2^{n_2}, x_3^{n_3})$, and this is
a contradiction. Therefore $P\in T_3$, and the proof of the claim
is complete.

\end{proof}

Let $T_4$ be the set of the union of $T_3$ and the set of all
minimal prime ideal $(x_1, x_2, x_3, x_4)$. It is clear that $T_4$
is a finite set. By [HH2, Lemma 3.2], we can choose an positive
integer $n$ such that for every prime ideal $P\in T_4$
$$((x_1^{n_1}, x_2^{n_2}, x_3^{n_3}): x_4^{n_4}))_P\subseteq ((x_1^{n_1}, x_2^{n_2}, x_3^{n_3}):
x_j^n))_P.$$ for all positive integers $n_1, n_2, n_3, n_4$ and
$j\leq 3$. Replacing $x_j$ by $x_j^n$, we have
$$((x_1^{n_1}, x_2^{n_2}, x_3^{n_3}): x_4^{n_4}))_P\subseteq ((x_1^{n_1}, x_2^{n_2}, x_3^{n_3}):
x_j))_P. \eqno (4.4)$$ for all positive integers $n_1, n_2, n_3,
n_4$ and $j\leq 3$.

For $P\notin T_4$, if $Pc\subset (x_1^{n_1}, x_2^{n_2},
x_3^{n_3}): x_j)$ for some element $c\in R$, then $Px_jc\in
(x_1^{n_1}, x_2^{n_2}, x_3^{n_3})$. Thus $x_jc\in (x_1^{n_1},
x_2^{n_2}, x_3^{n_3})$ by the choice of $P$. This shows that, for
all positive integers $n_1, n_2, n_3$, the associated primes of
$R/(x_1^{n_1}, x_2^{n_2}, x_3^{n_3}): x_j)$ lie in $T_4$. From
this fact and (4.4), one can conclude easily that
$$((x_1^{n_1}, x_2^{n_2}, x_3^{n_3}): x_4^{n_4})\subseteq ((x_1^{n_1}, x_2^{n_2}, x_3^{n_3}):
x_j). \eqno (4.5)$$
 hold for all positive integers $n_1, n_2,
n_3, n_4$ and $j\leq 3$.

\begin{claim}{4.}
For every permutation $i_1, i_2, i_3$ of $1, 2, 3$, we have

$$((x_{i_1}^{n_1}, x_{i_2}^{n_2}, x_4^{n_4}): x_{i_3}^{2})=
((x_{i_1}^{n_1}, x_{i_2}^{n_2}, x_4^{n_4}): x_{i_3}^{3}). $$ hold
 for all positive integers $n_1, n_2, n_4$.
\end{claim}
\begin{proof}
We only prove the following case:
$$((x_1^{n_1}, x_2^{n_2}, x_4^{n_4}): x_3^2)\subseteq
((x_1^{n_1}, x_2^{n_2}, x_4^{n_4}): x_3^{3}).$$

In fact, for any element $c\in ((x_1^{n_1}, x_2^{n_2}, x_4^{n_4}):
x_3^{3})$, we may express $$x_3^{3}c= x_1^{n_1}c_1+ x_2^{n_2}c_2+
x_4^{n_4}c_4$$ for some elements $c_1, c_2, c_4\in R$. So by
(4.5), there exists $c_3\in R$ such that $
x_3^{3}(x_3c-x_4^{n_4}c_3)\in (x_1^{n_1}, x_2^{n_2})$. By (4.3),
we have $x_3^2c\in (x_1^{n_1}, x_2^{n_2}, x_4^{n_4})$, and this
proves the claim.
\end{proof}

Replacing $x_1, x_2, x_3$ by $x_1^2, x^2_2, x_3^2$ if necessary,
we may assume
$$((x_{i_1}^{n_1}, x_{i_2}^{n_2}, x_4^{n_4}): x_{i_3})=
((x_{i_1}^{n_1}, x_{i_2}^{n_2}, x_4^{n_4}): x_{i_3}^{2}) \eqno
(4.6)$$ for every permutation $i_1, i_2, i_3$ of $1, 2, 3$ and for
all positive integers $n_1, n_2, n_4$.

Moreover, we can replace $x_4$ by a power of $x_4$ if necessary,
and assume that $((x_1^2, x_2^2, x_3^2): x_4)=((x_1^2, x_2^2,
x_3^2): x_4^2)$. Since $R$ is a excellent ring, we can choose an
element $x_5$ which is not contained in any minimal prime ideal of
$(x_1, x_2, x_3, x_4)$ such that $(R/(x_1^2, x_2^2, x_3^2,
x_4^2))_P$ is a $CM$ ring for every prime ideal $P$ with
$x_5\notin P$. By the choice of $x_4$ and Proposition 4.2, we
conclude that $(R/(x_1^2, x_2^2, x_3^2, x_4^{n_4}))_P$ is a $CM$
ring for every prime ideal $P$ with $(x_1, x_2, x_3, x_4)\subseteq
P$, $x_5\notin P$ and all positive integers $n_4$. It follows from
(4.6) and Proposition 4.2 that $(R/(x_1^{n_1}, x_2^{n_2},
x_3^{n_3}, x_4^{n_4}))_P$ are $CM$ local rings for such prime
ideals $P$ and all positive integers $n_1, n_2, n_3, n_4$.

\begin{claim}{5.}
Let $m$ be a maximal ideal of $R $ such that $m$ does not contain
the ideal $(x_1, x_2, x_3, x_4, x_5)$. Then for every $x_j$
($1\leq j\leq 3$), $x_j^5\text{H}_m^i(R)=0$ for $i< \text{ht}(m)$.
\end{claim}

\begin{proof}
Note that if $m$  does not contain the ideal $(x_1, x_2, x_3,
x_4)$, then the conclusion of the claim follows from Claim 2. So
we may assume that $(x_1, x_2, x_3, x_4)\subseteq m$ and
$x_5\notin m$. Replacing $R$ by $R_m$, we  assume that $R$ is a
local ring with the maximal ideal $m$. Observe that $x_1, x_2$ is
a regular sequence, (4.3) and (4.5). We have for $1\leq i\leq 4$
$$x_j((x_1^{n_1}, x_2^{n_2}, \cdots, x_{i-1}^{n_{i-1}}) :
x_i^{n_i})\subset (x_1^{n_1}, x_2^{n_2}, \cdots,
x_{i-1}^{n_{i-1}})$$ for each j with $1\leq j\leq 3$. Moreover, by
the choice of $m$,  $R/(x_1^{n_1}, x_2^{n_2}, x_3^{n_3},
x_4^{n_4})$ are $CM$ local rings for  all positive integers $n_1,
n_2, n_3, n_4$. So the conclusion of the claim follows from Lemma
4.5, and this ends of the proof of the claim.
\end{proof}

Now,  if  $\text {dim} (R/(x_1, x_2, x_3, x_4,
 x_5))=0$, then  there are only
a finite number of maximal ideals $m$ with  $(x_1, x_2, x_3, x_4,
x_5)\subseteq m$, So by Claim 5, and Corollary 4.4, a power of
$x_j$ ($1\leq j\leq 3$) is a uniform local cohomological
annihilator  of $R$. In particular, we have proved the theorem in
the case $d=5$.
\end{proof}

Before the end of the paper, we give a remark on the technique
used in this section.
\begin{remark}
The proof of Theorem 4.6 depends heavily on    Goto's result
(Proposition 4.2). The condition $(0 : x)=(0 : x^2)$ in
Proposition 4.2 is very restricted if one considers a lot of local
rings at the same time. We explain this more explicitly by means
of the proof of Theorem 4.6. Let  $R$, $x_1, x_2, x_3, x_4$ and
$x_5$ be as chosen as in the proof of Theorem 4.6. Although we can
find an element $x_6$ such that for all maximal ideals $m$ with
$x_6\notin m$, $R/(x_1^{n_1}, x_2^{n_2}, x_3^{n_3}, x_4^{n_4},
x_5)$ are $CM$ rings for all positive integers $n_1, n_2, n_3,
n_4$, it is very difficult to check that

$$((x_1^{n_1}, x_2^{n_2}, x_3^{n_3},
x_4^{n_4}): x_5)= ((x_1^{n_1}, x_2^{n_2}, x_3^{n_3}, x_4^{n_4}):
x_5^2)$$ hold for all positive integers $n_1, n_2, n_3, n_4$. So
one can not use Proposition 4.2 to conclude that $R/(x_1^{n_1},
x_2^{n_2}, x_3^{n_3}, x_4^{n_4}, x_5^{n_5})$ are $CM$ rings for
all positive integers $n_1, n_2, n_3, n_4, n_5$. Thus the method
of this paper can not be used to solve the remaining case of
Conjecture 4.1.  However, with a little more effort, we can prove
that, after replacing  $x_1, x_2, x_3, x_4, x_5$ by suitable
powers of them,   for each j ($1\leq j\leq 3$)
$$x_j((x_1^{n_1}, x_2^{n_2}, x_3^{n_3},
x_4^{n_4}): x_5^{n_5})\subseteq (x_1^{n_1}, x_2^{n_2}, x_3^{n_3},
x_4^{n_4})$$ hold for all positive integers $n_1, n_2, n_3, n_4,
n_5$. Due to this fact, one can prove easily that for any
excellent domain $R$ and an element of $x$ in $R$ such that $R_x$
is $CM$,  there exists an positive integer $n$,
$x^n\text{H}_m^i(R)=0$ for every maximal ideal $m$ and $i< \text
{min}(5, \text{ht}(m))-1$.
\end{remark}

{\bf Acknowledgment:}\hspace{0.2cm} The author is deeply grateful
to the referee for his or her useful pointed comments on the paper
and for drawing the author' attention to the results about the
annihilation of local cohomology modules.

\vspace* {0.3cm}

\begin{list}
\bf {\textbf{References}}{} \small \item[\rm{[Ba]}]  Bass, H.:  On
the ubiquity of Gorenstein rings.
 {\it Math. Z.}  {\bf 82},  8-28 (1963).

\item[\rm{[BE]}]   Buchsbaum D., Eisenbud, D.: What makes a
complexes exact. {\it J. Algebra}  {\bf 25}, 259-268(1973).

\item[\rm{[BRS]}]  Brodmann, M., Rotthaus, C., Sharp, R. Y.:  On
annihilators and associated primes of local cohomology modules.
 {\it J. Pure Appl. Algebra}  {\bf 153},  197-227 (2000).

\item[\rm{[BS]}]  Brodmann, M.,  Sharp, R. Y.:   Local cohomology:
an algebraic introduction with geometric applications  .
 {\it Cambridge Univ. Press, Cambridge},(1998).

\item[\rm{[Fa1]}]  Faltings, G.:  Uber die Annulatoren
Kohomologiegruppen. {\it Arch. Math.} {\bf 30}, 473-476(1978).

\item[\rm{[Fa2]}]  Faltings, G.:  Der Endlichkeitssatz in der
lokalen Kohomologie. {\it Math. Ann.} {\bf 255}, 45-56(1981).

\item[\rm{[Go]}]   Goto, S.:  Approximately Cohen-Macaulay rings.
{\it J. Algebra} {\bf 76}, 214-225(1982).

\item[\rm{[Gr]}]  Grothendieck, A.:  Local Cohomology.
 {\it Lect. Notes Math. } {\bf 41},  Springer-Verlage, New York
 (1963).

\item[\rm{[HH1]}]  Hochster, M., Huneke C.:  Tight closure,
invariant theory, and the Briancon-Skoda theorem. {\it J. Amer.
Math. Soc.} {\bf 3}, 31-116(1990).

\item[\rm{[HH2]}]  Hochster, M., Huneke C.:  Infinite integral
extensions and big Cohen-Macaulay algebras {\it Annals of math.. }
{\bf 3}, 53-89(1992).

\item[\rm{[HH3]}]  Hochster, M., Huneke C.: Applications of the
existence of big  Cohen-Macaulay algebras. {\it Adv. Math. } {\bf
113}, 45-117(1995).

\item[\rm{[HH4]}] Hochster, M., Huneke C.: Comparison of symbolic
 and ordinary powers of ideals. {\it Invent. Math. } {\bf 147},
349-369(2002).

\item[\rm{[HK]}]   Herzog, J. Kunz,  E.: Der Kanonische Module
eines Cohen-Macaulay Rings. {\it Lecture Notes in Mathematics}
 {\bf 238}  Springer-Verlage, New York(1971).

\item[\rm{[Hu]}]Huneke, C.: Uniform bounds in noetherian rings.
{\it Invent. Math. } {\bf 107}, 203-223(1992).

\item[\rm{[KS]}]  Khashyarmanesh, K.  and Salarian, S.: Faltings
theorem for  the  annihilation of local cohomology modules over a
Gorestein ring, {\it Proc. AMS  } {\bf 132}, 2215-2220(2004).

\item[\rm{[Ma]}] Matsumura, H.: Commutative Ring Theory, {\it
Cambrige Univ. Press, New York},(1986).

\item[\rm{[Ra1]}]  Ratliff, L.:  On quasi-unmixed local domain,
the altitude formula , and the chain condition for prime ideals I.
{\it Amer. J. Math.} {\bf 91}, 508-528(1976) 508-528(1969).

\item[\rm{[Ra2]}]  Ratliff, L.:  On quasi-unmixed local domain,
the altitude formula , and the chain condition for prime ideals
II. {\it Amer. J. Math.} {\bf 92}, 99-144(1970).

 \item[\rm{[Ra3]}]
Ratliff, L.:  Characterizations of catenary rings {\it Amer. J.
Math.} {\bf 93}, 1070-1108(1971).

\item[\rm{[Ra4]}]  Ratliff, L.: Catenary rings and the altitude
formula {\it Amer. J. Math.} {\bf 92}, 458-466(1972).

\item[\rm{[Rag1]}]  Raghavan K.:   Uniform annihilation of local
cohomology and of Koszul homology, {\it Math. Proc. Camb. Phil.
Soc.} {\bf 112}, 487-494(1992).

\item[\rm{[Rag2]}]  Raghavan K.:   Local-global principle for
annihilation of local cohomology, {\it Contemporary  Math.} {\bf
159}, 329-331(1994).

\item[\rm{[Ro]}]  Roberts, P.:  Two applications of dualizing
complexes over local rings. {\it Ann. Sci. $\acute{E}$c. Norm.
Sup$\acute{e}$r.}  {\bf 9}, 103-106(1976).

\item[\rm{[Rot]}] Rotman, J.: An introduction to homological
algebra, {\it Pure Appl. Math.,
 Academic Press,  New York},(1979).

\item[\rm{[Sc]}]  Schenzel, P.:   Cohomological annihilators. {\it
Math. Proc. Camb. Phil. Soc.}  {\bf 91}, 345-350(1982).

\item[\rm{[Sh1]}]  Sharp, R.:  Local cohomology theory in
commutative algebra. {\it Q. J.  Math. Oxford}   {\bf 21},
425-434(1970).

\item[\rm{[Sh2]}]  Sharp, R.:  Dualizing complexes for commutative
noetherian rings. {\it   Math. Proc. Camb. Phil. Soc.}  {\bf 78},
369-386(1975).

\end{list}

\end{document}